\documentclass{aptpub}

\newcommand{\R}{\mathbb{R}}
\newcommand{\N}{\mathbb{N}}
\newcommand{\Z}{\mathbb{Z}}

\newcommand{\Pp}{\mathbb{P}}

\newcommand{\Ee}{\mathcal{E}}

\usepackage{amsmath}
\usepackage{esvect}

\usepackage{color}
\usepackage{comment}

\DeclareMathOperator{\Sol}{Sol}

\usepackage[colorlinks = true,linkcolor = blue,urlcolor  = blue,citecolor = blue,anchorcolor = blue]{hyperref}

\authornames{PABLO A. GOMES, ALAN PEREIRA AND REMY SANCHIS} 
\shorttitle{Upper bounds for critical probabilities in Bernoulli Percolation models} 

\begin{document}

\title{Upper bounds for critical probabilities in \\ Bernoulli Percolation models} 

\authorone[Universidade de São Paulo]{Pablo A. Gomes} 
\authortwo[Universidade Federal de Alagoas]{Alan Pereira} 
\authorthree[Universidade Federal de Minas Gerais]{Remy Sanchis}

\addressone{Rua do Matao, 1010 - Cidade Universitária, São Paulo - SP, Brazil. 05508-090} 
\emailone{pagomes@usp.br} 
\addresstwo{BR-104, Km 97 - Cidade Universitária, Maceió - AL, Brazil. 57072-970} 
\emailtwo{alan.pereira@im.ufal.br} 
\addressthree{ Av. Pres. Antônio Carlos, 6627, Belo Horizonte - MG, Brazil. 31270-901} 
\emailthree{rsanchis@mat.ufmg.br} 

\begin{abstract}
In this paper we study bond and site Bernoulli Percolation models on $\Z^d$ for $d \geq 3$ with parameter $p$, both in the oriented and non-oriented versions. The main macroscopic quantity of interest is the probability of long-range order and the existence of a non-trivial threshold is well established. Precise numerical results for the threshold values are available in the literature, but mathematically rigorous bound are mostly  restricted to 2D lattices. Utilizing dynamical coupling techniques, we introduce a comprehensive set of new rigorous upper bounds that corroborate existing numerical values.
\end{abstract}



\section{Introduction} The study of Bernoulli  percolation on $\Z^d$ is more than 60 years old and the existence of a non-trivial phase transition for $d\geq 2$ is well established for the model and several of its variants, but the exact value of the critical parameter $p_c$ is seldom known. A celebrated result of Kesten (see \cite{K3}) proved that the critical probability in Bernoulli bond percolation on $\Z^2$ is $1/2$. Besides that, the exact critical probability was established only for some planar lattices where symmetry and duality property are well explored.  
On the other hand, for dimensions $d\geq 3$, there is not much hope of finding exact values for the critical probability. Several non-rigorous estimates, either via Monte Carlo Methods (see the very efficient algorithm in \cite{NZ}), or statistical estimates based on a comparison with dependent percolation (see Section 6.2 of \cite{BR} for an overview), are present in the literature, but rigorous bounds are rare (see, in different contexts, \cite {CR, L,W,WP}). In this paper, we focus on finding rigorous upper bounds for site and bond Bernoulli Percolation on $\Z^d$ for every $d\geq 3$, in both the oriented and the non-oriented cases. Although those bounds are still far from the values obtained from Monte Carlo simulations, we believe that most of them are the best rigorous upper bounds in the literature.

The primary tools we use are dynamical couplings between the models we seek to understand and models where bounds or precise values for the critical probabilities are known.  Such couplings are now a classic tool in percolation theory, but its extensive use to establish upper bounds on critical parameters for all dimensions, we believe, is new. In \cite{CS}, the authors establish inequalities between the critical thresholds of site and bond non-oriented percolation through the use of differential inequalities.

The remainder of the text is organized as follows: in Section \ref{sec:model} we  precisely define the models; state the main results; and, in Subsection \ref{sec:comments} we give a numerical table with upper bounds for the critical probabilities in homogeneous Bernoulli percolation models for dimensions up to $d=9$. In Section \ref{sec:couplings} we establish some dynamical couplings and in Section \ref{sec:proofs} we prove the theorems.

\section{The models and main results}\label{sec:model}

For $d \geq 1$, the underlying graph for all the models will be the $d$-dimensional hypercubic lattice in which the set of vertices is $\Z^d$ and the set of edges is the set of non ordered pairs $E(\Z^d) : = \{ \langle v, u \rangle : v, u \in \Z^d ~ \text{and} ~ |v-u| =1\}$. We abuse notation and denote this graph simply by $\Z^d$. Throughout the text, for each $x \in \R$, we will denote by $\lfloor x \rfloor$ the largest integer smaller or equal to $x$.

This section is divided into two subsections. The first is devoted to bond percolation models, while the second is devoted to site percolation models.

Before we proceed with formal definitions of the models, we introduce a notation that will be useful during the text.

\begin{definition}\label{Rsol}
Given two real numbers $a \in (0,1)$ and $b > 1$, we denote by $\Sol(a,b)$, the unique solution of $x = a[x + (1-x)^b]$ in the interval $(0,1)$.
\end{definition}

     Observe that for every $x$ in the interval $(0,1)$, the right hand side of the equation is smaller than $a$, therefore $\Sol(a,b) < a$. Denoting $f(x) = x - a[x + (1-x)^b]$, we have that $f(0) = -a < 0$ and $f(1) = 1 - a > 0$. Moreover, observe that $f'(x) = (1 - a) + ab(1-x)^{b-1} > 0$, for every $x \in (0,1)$, what guarantees the uniqueness of $\Sol(a,b)$.

\medskip

\subsection{Bond percolation}

Given  $p \in [0,1]$, consider a family of independent random variables $\{X_e\}_{e \in E(\Z^d)}$, where, for each $e \in E(\Z^d)$, $X_e$ has Ber($p$) distribution. Let $\mu_e$ be the law of $X_e$, and let $\Pp_p := \prod_{e \in E(\Z^d)} \mu_e$ be the resulting product measure. We declare an edge $e$ to be {\it open} if $X_e=1$ and {\it closed} otherwise. 

We first consider the non-oriented case and denote by $\{x \leftrightarrow y\}$ the event where  $x, y \in \Z^d$ are connected by an open path, i.e., there exist $x_0, \dots, x_n$ such that $x_0 = x$, $x_n =y$ and each $\langle x_{j-1}, x_j \rangle$ belongs to $ E(\Z^d)$ and is open for $j=1, \dots, n$. Let  $\mathcal{C}^b_0 := \{ x \in \Z^d:  0 \leftrightarrow x\}$ be the open cluster of the origin, and $|\mathcal{C}^b_0|$ its size.

We define the percolation probability by $\theta^{b}_d(p) := \Pp_p(|\mathcal{C}^b_0|= \infty)$. The critical point for the {\it non-oriented bond} Bernoulli percolation model will be denoted by
 \begin{equation*} \label{eq:defCriticalPoint}
    p_c^{b}(d)=\sup\{p \geq 0 : \theta^b_d(p)=0\}.
\end{equation*}

We now consider the oriented case. Let $\{e_1, \dots, e_d\}$ be the set of positive unit vectors of $\Z^d$. We denote by $\{x \rightarrow y\}$ the event where  $x, y \in \Z^d$ are connected by an oriented open path, i.e., there exist $x_0, \dots, x_n$ such that $x_0 = x$, $x_n =y$ and for each $j=1, \dots, n$, we have $x_j - x_{j-1} \in \{e_1, \dots, e_d\}$ and $\langle x_{j-1}, x_j\rangle$ is open. Let $\vv{\mathcal{C}}^b_0 := \{ x \in \Z^d : 0 \rightarrow x\}$ be the oriented open cluster of the origin, and $|\vv{\mathcal{C}}^b_0|$ its size.

Analogously, we define 
$\vv{\theta}^{\, b}_d(p) :=  \Pp_p(|\vv{\mathcal{C}}^b_0|= \infty)$ the corresponding oriented percolation probability, and we denote the critical point for the {\it oriented bond} Bernoulli percolation model by \[\vv{p}_c^{\, b}(d)=\sup\{p \geq 0 : \vv{\theta}^{\,b}_d(p)=0\}.\]

For non-oriented bond percolation model, our result is the following:

\begin{thm}\label{theo:bounds} Consider non-oriented bond Bernoulli percolation and 
let $p^{\ast}(d)$ be the unique solution in $(0,1)$ of 
\begin{equation} \label{eq:Theorem1}
  \prod_{i=0}^2 \left(1-(1-p)^{\lfloor\frac{d+i}{3}\rfloor}\right)=2-\sum_{i=0}^2(1-p)^{\lfloor\frac{d+i}{3}\rfloor}.  
\end{equation}

Then, for every $d\geq 3$, we have $p_c^b(d) \leq p^{\ast}(d)$.
\end{thm}

To state our result for oriented bond percolation model, we recall Definition~\ref{Rsol}. 

\begin{thm}\label{theo:bondO}
Consider oriented bond Bernoulli percolation on $\Z^d$. The three following statements hold.

1) If $d$ is even, then $ \vv{p}^{\,b}_c(d)\leq 1-(1/3)^{2/d}$;

2) For $d \geq 4$, we have that \[\vv{p}^{\,b}_c(d)\leq \frac{1}{d}+\frac{C_d}{d^2},\]
where \[C_d=1+\frac{8}{d}+\frac{d^{5/2}}{(\sqrt{2\pi})^{d-1}}\left[\frac{d-1}{d-3}\right]e^{\frac{1}{12d}}.\]

3) For any dimension $d\geq 2$, we have that  $\vv{p}^{\,b}_c(d+1)\leq \Sol\big(\vv{p}^{\,b}_c(d), (d+1)/d \big)$.

\end{thm}

\begin{rem} The second item of Theorem~\ref{theo:bondO} is implicitly given in \cite{GPS}. A precise description is given in the proof.
\end{rem}

\begin{rem} It is known (see \cite{CD}) that $\vv{p}^{\,b}_c(d)\sim 1/d$, hence the second upper bound above is asymptotically sharp.
\end{rem}

\subsection{Site percolation}

Given a parameter $p \in [0,1]$, we consider a family $\{X_v\}_{v \in \Z^d}$ of independent Bernoulli random variables with parameter $p$. As before, $\Pp_p$ will denote the resulting product measure. A vertex $v \in \Z^d$ is declared to be {\it open} if  $X_v = 1$ and {\it closed} otherwise. 

Now, the sequence of vertices $(x_0, \dots, x_n)$ is said to be an open path if all the vertices are open and for each $j = 1, \dots, n$, $\langle x_{j-1}, x_j\rangle \in E(\Z^d)$. If in addition to these conditions, for each $j = 1, \dots, n$, we have $x_j - x_{j-1}\in \{e_1, \dots, e_d\}$, the sequence is said to be an oriented open path. Having these definitions we  define $\mathcal{C}^s_0$, $\vv{\mathcal{C}}^s_0$, $\theta_d^{\, s}(p)$ and $\vv{\theta}_d^{\, s}(p)$ accordingly.

Finally, the critical points for {\it non-oriented} and {\it oriented} site Bernoulli percolation are respectively given by
\[ p_c^{s}(d)=\sup\{p \geq 0 : \theta^s_d(p)=0\} \quad \quad \text{and} \quad \quad \vv{p}_c^{\, s}(d)=\sup\{p \geq 0 : \vv{\theta}^{\,s}_d(p)=0\} .\]

To state our results for site percolation models, we recall Definition~\ref{Rsol}. 


\begin{thm}\label{theo:siteNO}
Consider non-oriented site Bernoulli percolation on $\Z^d$. The three following statements holds.

1) If $d$ is even, then $ p_c^{s}(d)\leq 1-(0,32)^{2/d}$; 

2) If $d$ is divisible by 3, then $p_c^{s}(d)\leq 1-(1/2)^{3/d}$;

3) For any dimension $d\geq 2$, we have that $p_c^{s}(d+1)\leq \Sol(p^{s}_c(d),2d/(2d-1))$.


\end{thm}

\begin{thm}\label{theo:siteO}
Consider oriented site Bernoulli percolation on $\Z^d$. The two following statements holds.

1) If $d$ is even, then $ \vv{p}_c^{\,s}(d) \leq 1-(1/4)^{2/d}$;

2) For any dimension $d \geq 2$, we have that $\vv{p}_c^{\,s}(d+1) \leq \Sol\big(\vv{p}^{\,s}_c(d),(d+1)/d \big)$. 


\end{thm}

\subsection{Explicit bounds}\label{sec:comments}

In this subsection, we present a table with rigorous upper bounds for bond and site  Bernoulli Percolation in $\Z^ d$, up to $d=9$, in both the oriented and the non-oriented cases. 

In each column we give a numerical upper bound rounded to four decimals for the critical probability in the head of the column.

\begin{table}[h!]
    \centering

\bigskip
   
\begin{tabular}{|c|c|c|c|c|}  

\hline
\textbf{Dimension} & $p_c^b(d)$ & $\vv{p_c}^b(d)$ & $p_c^s(d)$ &   $\vv{p_c}^s(d)$\\ \hline

3 & 0,3473 $^{(1)}$ &0,5680 $^{(2)}$&0,5000 $^{(*)}$&0,6422 $^{(9)}$\\   \hline

4 & 0,2788 $^{(1)}$ &0,4227 $^{(3)}$&0,4344 $^{(6)}$&0,5000 $^{(10)}$ \\   \hline

5 & 0,2284 $^{(1)}$ & 0,3926 $^{(4)}$&0,4156 $^{(8)}$& 0,4615 $^{(11)}$\\   \hline

6 & 0,1922 $^{(1)}$& 0,2734 $^{(5)}$ & 0,2929 $^{(7)}$& 0,3701 $^{(10)}$\\   \hline

7 & 0,1682  $^{(1)}$& 0,2028 $^{(5)}$& 0,2866 $^{(8)}$& 0,3533 $^{(11)}$\\   \hline

8 & 0,1486 $^{(1)}$&0,1627  $^{(5)}$& 0,2479 $^{(6)}$ &0,2929 $^{(10)}$\\   \hline

9 & 0,1326 $^{(1)}$ &0,1371 $^{(5)}$&0,2063 $^{(7)}$&0,2844 $^{(11)}$\\   \hline

\end{tabular}

\caption{Numerical bounds up to dimension 9.  .}
   
\label{tab:my_label}
  
\end{table}


\vspace{0.5cm}

To obtain each bound, we used the following:

\begin{enumerate}
\item[$(\ast)$] The bound for non-oriented site percolation in dimension $d = 3$, $p_c^s(3) < 1/2$,  was given in \cite{CR}
    \item[(1)]  All bounds for non-oriented bond percolation were obtained computing $p^{\ast}(d)$ as given in Theorem~\ref{theo:bounds}.
    \item[(2)] The bound for oriented bond percolation in dimension $d = 3$ was obtained from Item~3) of Theorem~\ref{theo:bondO} together the bound $\vv{p}^{\,b}_c(2)\leq 2/3$ given in \cite{L}.
    \item[(3)] The bound for oriented bond percolation in dimension $d = 4$ was obtained from Item~1) of Theorem~\ref{theo:bondO}.
    \item[(4)] The bound for oriented bond percolation in dimension $d = 5$ was obtained from Item~3) of Theorem~\ref{theo:bondO} together the bound obtained for $d=4$.
    \item[(5)] The bounds for oriented bond percolation in dimensions $d \in \{6, 7, 8, 9\}$ were obtained from Item~2) of Theorem~\ref{theo:bondO}. 
    \item[(6)] The bounds for non-oriented site percolation in dimensions $d \in \{4, 8\}$ were obtained from Item~1) of Theorem~\ref{theo:siteNO}. 
    \item[(7)] The bounds for non-oriented site percolation in dimensions $d \in \{6, 9\}$ were obtained from Item~2) of Theorem~\ref{theo:siteNO}. 
    \item[(8)] The bound for non-oriented site percolation in dimensions $d \in \{5, 7\}$ were obtained from Item~3) of Theorem~\ref{theo:siteNO} together the bounds obtained for $d \in \{4,6\}$.
    \item[(9)]  The bound for oriented site percolation in dimension $d = 3$ was obtained from Item~3) of Theorem~\ref{theo:bondO} together the bound $\vv{p}^{\,s}_c(2)\leq 3/4$ given in \cite{L}.
     \item[(10)] The bounds for oriented site percolation in dimensions $d \in \{4,6,8\} $ were obtained from Item~1) of Theorem~\ref{theo:siteO}. 
    \item[(11)]The bounds for oriented site percolation in dimensions $d \in \{5,7,9\}$ were obtained from Item~2) in Theorem~\ref{theo:siteO} together the bounds obtained for $d \in \{4,6,8\}$ 
\end{enumerate}

\begin{rem} Recently, Yu and Wierman, using similar methods, published an article (see \cite{YW}) proving the same upper bound for non-oriented bond percolation in $\Z^3$ as ours.
\end{rem}

\section{The Dynamical Couplings}\label{sec:couplings}

Although we are only interested in homogeneous percolation, some of the tools we will use are related to couplings between {\it anisotropic} bond percolation models. The algorithm used for establishing the couplings is reminiscent of the one used in Lemma~1 of \cite{GM}.

We now define the anisotropic (or inhomogeneous) non-oriented and oriented bond percolation models. For each $i=1, \dots, d$, let $E_i = \{ \langle x, x + e_i \rangle: x \in \Z^d\}$ be the set of edges parallel to $e_i$. Given $p_1, \dots, p_d \in [0,1]$,  we consider a family of independent random variables $\{X_e\}_{e \in E(\Z^d)}$, but now, for each $e \in E_i$, $X_e$ has Ber($p_i$) distribution, $i=1, \dots, n$. The open cluster of the origin $\mathcal{C}^b_0$ and the oriented open cluster of the origin $\vv{\mathcal{C}}^b_0$ are  defined analogously. The probabilities that $|\mathcal{C}^b_0|$, $|\vv{\mathcal{C}}^b_0|$ are infinite, will be denoted by $\theta_d(p_1, \dots, p_d)$ and $\vv{\theta}_d (p_1, \dots, p_d)$ respectively.

The first coupling is the content of Proposition 1 in \cite{GPS2}:
\begin{prop} \cite[Proposition 1]{GPS2} \label{prop:coupling} 
Consider inhomogeneous non-oriented Bernoulli bond percolation in $\Z^d$. Let $p_1, \dots, p_{d+1} \in [0,1]$ and let $\Tilde{p}_d \in [0,1]$ be such that \[(1-\Tilde{p}_d) = (1-p_d)(1-p_{d+1}).\] Then,
    $\theta_{d+1}(p_1, \dots,p_d,  p_{d+1})  \geq \theta_{d}(p_1, \dots, p_{d-1}, \Tilde{p}_{d}) .$
\end{prop}

Observe that throughout $d-3$ application of Proposition~\ref{prop:coupling} we are able to compare a $d$-dimensional model with a three dimensional one. In the next proposition we will obtain a comparison of the model in $\Z^3$ with the model in the triangular lattice $\mathbf{T}$ (defined below). Since the critical surface for inhomogeneous bond percolation in $\mathbf{T}$ is well established (see Theorem 11.116 in \cite{Grim}), this comparisons will be our strategy to obtain a proof of Theorem~\ref{theo:bounds}.

Now we define the triangular lattice $\mathbf{T}$. This lattice is simply $\Z^2$ with  extra edges of the form $\langle v,v+(1,1)\rangle$. That is, $\mathbf{T} = (V_{\mathbf{T}}, E_{\mathbf{T}})$ where the set of vertices is given by $V_{\mathbf{T}} = \Z^2$ and the set of edges $E_{\mathbf{T}}$ is the set of non ordered pairs $\{\langle v, u\rangle : v-u = (1,0), (0,1) \text{ or } (1,1) \}$. We will denote by $\theta_{\mathbf{T}}$  the corresponding percolation probability. 

We will consider inhomogeneous bond percolation on $\mathbf{T}$, where the corresponding parameters $(p_1,p_2,p_3)$ will refer to edges of the form $\langle v,v+(1,0)\rangle$, $\langle v,v+(0,1)\rangle$ and $\langle v,v+(1,1)\rangle$, respectively.

In the next proposition, we construct a monotonic coupling between inhomogeneous bond percolation on the  triangular lattice $\mathbf{T}$ with parameters $(p_1,p_2,p_3)$ and on $\Z^3$ with the same parameters.

Our proof is reminiscent of the  arguments used in~\cite{CR}, where a strictly inequality is given for the homogeneous site model. In the homogeneous setting, the next proposition is just a special case of Theorem 1 in~\cite{BS}. We observe that this theorem already gives the bound $p^b_c(3)\leq p^b_c(\mathbf{T})\leq 0.3473 $. Moreover, our proof is essentially the same as in \cite{BS} and we include it here for sake of completeness since we need to deal with the inhomogeneous setting. More general results involving strict inequalities for the critical parameter of quotient graphs have been recently obtained in \cite{MS}.

\begin{prop}\label{prop:triangular}
Given the parameters $(p_1,p_2,p_3)\in [0,1]$, consider two inhomogeneous bond Bernoulli percolation processes: on the triangular lattice  and on $\Z^3$. Then
\[\theta^b_{\mathbf{T}}(p_1,p_2,p_3)\leq\theta_3(p_1,p_2,p_3).\]
\end{prop}

\begin{proof}
We will construct a dynamic coupling between the percolation process on $\Z^{3}$ with parameters $(p_1,p_2, p_3)$  and an infection process over $\mathbf{T}$.  We will do it in such a way that the law of infected sites in $\mathbf{T}$ is the same as the law of the open cluster of the origin for anisotropic percolation on $\Z^3$ and also that, if the infection process survives, the open cluster of the origin of the process in $\Z^3$ must be infinite. 

Before introduce the coupling, let us set the following auxiliary notation
\begin{equation}\label{eq:tau}
\tau(\pm(1,0))=\pm (1,0,0), \: \tau(\pm(0,1))=\pm (0,1,0) \mbox{ and } \tau(\pm(1,1))=\pm (0,0,1) .
\end{equation}

The coupling will be built based on a susceptible-infected strategy described as follows. First, we declare  the origin of $\mathbf{T}$ as the \textit{initial} infected component. Next, at each time-step, we possibly grow the infected component. 
Consider a vertex $v$ in the infected component of $\mathbf{T}$ and a neighbor $v+u$ out of the infected component. The vertex $v$ will be associated with some vertex $x(v)$ in the  open cluster of the origin in $\Z^3$. If $\langle x(v), x(v)+\tau(u) \rangle$,  is open, we infect $v+u$  (and write $x(v+u) = x(v) + \tau(u)$).


Formally, the coupling is obtained by the sequence of sets $(I_n, x(I_n), R_n, S_n)_{n \geq 0}$. Here, $I_n$ represents the \textit{infected vertices} in $\mathbf{T}$, $x(I_n)$ represents the vertices in $\Z^{3}$ associated with the infected vertices, and $R_n$ represents the \textit{removed edges} of $\mathbf{T}$. Finally, given $I_n$, $x(I_n)$ and $R_n$, the \textit{susceptible edges set} is given by
\begin{equation*}
    S_{n} := \{ \langle v, u \rangle \in E_{\mathbf{T}}: v \in I_{n} \: \mbox{and} \: u \notin I_{n} \} \cap R_{n}^C.
\end{equation*}
At time $n=0$, we set
\begin{itemize}
    \item $I_0 = \{0\} \subset V_{\mathbf{T}}$; 
    
    \item $R_0 = \emptyset \subset E_{\mathbf{T}}$;
    
    \item $x(0) = 0 \in \Z^{3}$;
    
    \item $S_0 := \{ \langle v, u \rangle : v \in I_{0} \: \mbox{and} \: u \notin I_{0} \} \cap R_{0}^C \:\: = \{ \langle 0, \pm (1,0)\rangle, \langle 0,\pm (1,0)\rangle ,\langle 0,\pm (1,1)\rangle \}$. 
\end{itemize}

This means that, at time $n=0$, only the vertex $0$ is infected, and it can potentially infect any of its neighbors, so all edges containing the origin are susceptible. After that, in each step, an infected vertex tries to infect a non-infected vertex through a susceptible edge (if the latter exists). From now on, we choose an arbitrary, but fixed, ordering of the edges in $\mathbf{T}$. Suppose that $I_n, x(I_n), R_n$ and $S_n$ are already defined. If there is no susceptible edge then the process stops. More specifically, if $S_n = \emptyset$, then for all $k \geq 1$,
    \begin{itemize}
        \item $I_{n+k} = I_n$;
        
        \item $R_{n+k} = R_n$;
        
        \item $S_{n+k} = S_n$.
    \end{itemize}
Otherwise, $S_n \neq \emptyset$. In this case, let $g_n$ be the smallest edge in $S_n$.
Since $g_n \in S_n$, it has to be equal to some $\langle v, v+u_n \rangle$, where $v \in I_n$, $v+u_n \notin I_n$, and $u_n \in \{ \pm (0,1), \pm (1,0), \pm (1,1)\}$. Then, $v$ \textit{infects} $v+u_n$ if $\langle x(v), x(v)+\tau (u_n) \rangle$ is open in $\Z^{3}$ (recall definition of $\tau$ in \eqref{eq:tau}). More precisely, if $\langle x(v), x(v)+\tau (u_n) \rangle$ is open in $\Z^{3}$ then we write
\begin{equation*}
    I_{n+1} := I_n \cup \{v+u_n\},
\end{equation*}
and define
\begin{equation*}
    x(v+u_n) := x(v) + \tau(u_n).
\end{equation*} 
Otherwise, if $\langle x(v), x(v)+ \tau(u_n) \rangle$ is closed in $\Z^{3}$,  we set $I_{n+1}:= I_n$.

In either case, we have explored $g_n$, thus we \textit{remove} it and write
\[ R_{n+1} := R_n \cup \{ g_n \}.\]
Next, to conclude our induction step, we set
\begin{equation*}
  S_{n+1} := \{ \langle v, u \rangle : v \in I_{n+1} \: \mbox{and} \: u \notin I_{n+1} \} \cap R_{n+1}^C.
\end{equation*}
Observe that the function $x: \cup_n I_n \longrightarrow \Z^{3}$ is injective.  
In fact, if $v= (v_1, v_2) \in I_n$, then by construction, we have that $x(v) = (x_1,x_2,x_3)$ satisfies
\[ v_1=x_1+x_3 \mbox{ and }
    v_2=x_2+x_3.\]
Now, observe that the image of $x$ is contained in the  open cluster of the origin $\mathcal{C}^b_0$ of $\Z^{3}$. Since $x$ is injective,  $|\cup_n I_n| \leq |\mathcal{C}^b_0|$.  Also, note that $\cup_{n} I_n$ has the same  law as  $\mathcal{C}^{\mathbf{T}}_0$, where $\mathcal{C}^{\mathbf{T}}_0$ is the  open cluster of the origin in $\mathbf{T}$ with parameters $p_1, p_2,p_3$.  Therefore,
\begin{equation*}
 \theta^b_{\mathbf{T}}(p_1,p_2,p_3) \leq \theta_{3}(p_1,p_2,p_3),
\end{equation*}
and the proof of Proposition \ref{prop:triangular} follows.
\end{proof}

The content of the next two propositions are devoted to construct two couplings that will be used to prove of Theorem~\ref{theo:bondO}, Theorem~\ref{theo:siteNO} and Theorem~\ref{theo:siteO}. The exploration methods in these couplings are reminiscent of the coupling in \cite{GSS} (where edges are split) and give an upper bound for the critical probability for any percolation model in $\Z^{d+1}$ as a function of the corresponding critical probability in $\Z^d$.

\begin{prop}\label{prop:crossover}
Consider oriented bond Bernoulli percolation on $\Z^{d+1}$ and suppose that 
\begin{equation} \label{eq:crossbond}
    p > \vv{p}^{\,b}_c(d) \left[p + (1-p)^{(d+1)/d}\right].
\end{equation}
Then, $\vv{\theta}^{\,b}_{d+1}(p)>0$.
\end{prop}


\begin{proof} 
We prove this proposition in two steps. First, we will define a multigraph $\Z^{d+1}_{\mathcal{E}}$ and show that, with certain parameters, the model on $\Z^{d+1}_{\Ee}$ is equivalent to the homogeneous model on $\Z^{d+1}$ with parameter $p$. To conclude, we will show that if $p$ satisfies Inequality \eqref{eq:crossbond}, then this new model on $\Z^{d+1}_{\Ee}$ dominates a supercritical model on $\Z^{d}$.

Let ${\Ee} = \{e_1, \dots, e_d\}$ denote the set of positive unit vectors of $\Z^d$. To avoid ambiguities, we denote the set of positive unit vectors of $\Z^{d+1}$ by $\{ u_1, \dots, u_{d+1}\}$. Consider the multigraph $\Z^{d+1}_{\Ee}$ defined as follows: the set of vertices is $\Z^{d+1}$ and the set of edges is given by $E(\Z^{d+1}_{\Ee}):= \left( \cup_{i=1}^d E_i \right) \cup E_{\Ee}$ where  $E_i = \{ \langle x, x + u_i \rangle: x \in \Z^{d+1}\}$, $i=1, \dots, d$, and we define $E_{\Ee} := \{\langle v, v+u_{d+1}\rangle_e : v \in \Z^{d+1} , e \in {\Ee}\}$. In words, each edge of $\Z^{d+1}$ parallel to $u_{d+1}$ is split into another $|{\Ee}|$ edges indexed by ${\Ee}$ in $\Z^{d+1}_{\Ee}$, while edges parallel to all other directions remain unmodified.

Consider now inhomogeneous oriented bond Bernoulli percolation on $\Z^{d+1}_{\Ee}$, where edges in $\cup_{i=1}^d E_i$ are open with probability $p$ and edges in $E_{\Ee}$ are open with probability $q$, where $q$ is such that $(1-p) = (1-q)^{|{\Ee}|}$. Clearly, the distribution of the  open cluster in this model is the same as in the homogeneous model on $\Z^{d+1}$ with parameter $p$.  We will construct a coupling showing that the model on $\Z_{\Ee}^{d+1}$ with parameters $p$ and $q$ as above, dominates the homogeneous model on $\Z^d$ with parameter $p/[1 - (1-p)q]$. 

First, for each $e_i \in \Ee$, let $ \sigma(e_i):= u_i \in \Z^{d+1}$. Then,  for each $v \in \Z^{d+1}$ and each $e \in {\Ee}$, let $A(v,e)$ be the event where either the edge $\langle v, v+\sigma(e) \rangle$ is open or, for some $k \geq 1$,
\begin{itemize}
    \item the edges $\langle v+ iu_{d+1}, v + (i+1)u_{d+1}  \rangle_e$, $i = 0, \dots, k-1$, are open and
    \item the edges $\langle v + iu_{d+1}, v + iu_{d+1} + \sigma(e) \rangle$, $i = 0, \dots, k-1$ are closed, and 
    \item $\langle v + ku_{d+1}, v + ku_{d+1} + \sigma(e) \rangle$ is open.
\end{itemize}
In the event  $A(v,e)$, we define $u(v,e) = v+\sigma(e)$ if $\langle v, v+\sigma(e) \rangle$ is open, or $u(v,e) = v + k u_{d+1} + \sigma(e)$ if $k \geq 1$ is such that the above three conditions are met. Observe  that
\begin{eqnarray*}
    \Pp(A(v,e)) &=& 
    p\sum_{i = 0}^{\infty} \left[q(1-p)\right]^i
    = \frac{p}{1 - (1-p)q}\\
    &=& \frac{p}{1-(1-p)(1-(1-p)^{1/d})} \\
    &=& \frac{p}{ p + (1-p)^{(d+1)/d} }.
\end{eqnarray*}

Similarly to what was done in Proposition \ref{prop:triangular}, we now build the sequence of sets $(I_n, x(I_n), R_n, S_n)_{n \geq 0}$. For $n =0$, we set
\begin{itemize}
    \item $I_0 = \{0\} \subset \Z^d$;
    \item $R_0 = \emptyset \subset E(\Z^d)$;
    \item $x(0) = 0 \in \Z_{\Ee}^{d+1}.$
\end{itemize}
Suppose that, for some $n \geq 0$, $I_n, x(I_n)$ and $R_n$ are already defined. Then $S_n \subset E(\Z^d)$ is given by
\begin{equation*}
    S_{n} := \{ \langle v, u \rangle : v \in I_{n} \: \mbox{and} \: u \notin I_{n} \} \cap R_{n}^C.
\end{equation*}
If $S_n = \emptyset$, our sequence becomes constant, that is, 
\begin{equation}\label{eq:stop}
    (I_{n+k}, x(I_{n+k}), R_{n+k}, S_{n+k}) = (I_n, x(I_n), R_n, S_n), \quad \forall k \geq 1.
\end{equation} 
Otherwise, let $g_n = \langle v , v + e \rangle$ be the smallest (according to a prefixed ordering, as in Proposition~\ref{prop:triangular}) edge of $S_n$, where $v \in I_n$, $e \in {\Ee}$, and $v+e \notin I_n$. We set $R_{n+1} = R_n \cup \{g_n\}$. We also set
\begin{equation*}
    I_{n+1} = 
    \begin{cases} I_n \cup \{v+e\}, &\text{if } A(x(v),e) \text{ occurs};\\
                  I_n, & \text{otherwise}.
    \end{cases}
\end{equation*}
In case  $A(x(v),e)$ occurs, we set $x(v+e) = u(x(v),e)$. 

Once our sequence $(I_n, x(I_n), R_n, S_n)_ {n \geq 0}$ is built, note that by construction, the function $x: \cup_n I_n \longrightarrow \Z^{d+1}_{\Ee}$ is injective. In fact, for each $v \in \cup_n I_n$, the projection of $x(v)$ into $\Z^d$ is equal to $v$. The conclusion of the proof follows in a similar way to that of Proposition \ref{prop:triangular}.
\end{proof}

\begin{prop}\label{prop:crossoversite}
Consider non-oriented site Bernoulli percolation on $\Z^{d+1}$ and suppose that 
\begin{equation}\label{eq:cross1}
p > p^s_c(d)  \left[p + (1-p)^{2d/(2d-1)}\right].
\end{equation}
Then $\theta^s_{d+1}(p)>0$.
\end{prop}


\begin{proof}
The strategy of the proof is similar to the previous proposition. As before, let $\{e_1, \dots, e_d\}$ denote the set of positive unit vectors of $\Z^d$ and let $\{ u_1, \dots, u_{d+1}\}$ denote the set of positive unit vectors of $\Z^{d+1}$. We will consider a graph where each vertex of $\Z^{d+1}$ will be split into another $2d-1$ vertices. 

For each $v \in \Z^{d+1}$, we define the set of split vertices $V_v := \{v^{(1)}, \dots, v^{(2d-1)}\}$. Let $\Z^{d+1}_V$ be the graph with vertex set  $\cup_{v \in \Z^{d+1}} V_v$ and  edge set
\begin{equation*}
    E(\Z^{d+1}_V) := \left\{ \langle x,y \rangle ~:~ x \in V_v \text{ and } y \in V_u, \text{ for some } \langle v,u \rangle \in E(\Z^{d+1}) \right\}.
\end{equation*}
Given $p$ satisfying \eqref{eq:cross1}, let $q$ be such that $(1-p) = (1-q)^{2d-1}$. We will consider the non-oriented site Bernoulli percolation model on $\Z^{d+1}_V$ with parameter $q$. Note that the probability that at least one vertex in $V_0$ percolates on $\mathbb{Z}^{d+1}_V$  equals  the probability that the origin percolates in the model of non-oriented site percolation on $\mathbb{Z}^{d+1}$ with parameter $p$.


To build the coupling between $\Z^{d+1}_V$ with parameter $q$ and $\Z^d$ with parameter $p/[1  - (1-p)q]$, we again construct a sequence of sets $(I_n, x(I_n), R_n, S_n)_{n \geq 0}$. Since we are  considering a site model, for each $n \geq 0$, the sets $R_n$ and $S_n$ will be, respectively, the sets of {\it removed} and {\it susceptible} vertices (instead of edges) at step $n$. We start with two infected vertices (otherwise the origin should have been split into $2d$ vertices, instead of $2d-1$). For $n = 0$, we set 
\begin{itemize}
    \item $I_0 = \{0, e_1\} \subset \Z^d$;
    \item $R_0 = \emptyset \subset \Z^d$;
    \item $x(0) = 0^{(1)} \in \Z_V^{d+1} \text{ and } x(e_1) = u_1^{(1)} \in \Z_V^{d+1}. $
\end{itemize}
Inductively, if for some $n \geq 0$, the sets $I_n, x(I_n)$ and $R_n$ are already defined, we set 
\begin{equation}\label{eq:SnSITE}
    S_n = \{u \in \Z^d : \langle v,u \rangle \in E(\Z^d) \text{ for some } v \in I_n\} \cap (I_n \cup R_n)^C.
\end{equation} 
If $S_n = \emptyset$, the sequence becomes constant, as in \eqref{eq:stop}. Otherwise, let $a_n$ be the smallest (in a preset ordering) vertex of $S_n$. We can write $a_n =  v+ e$, such that $v \in I_n$ and $e \in \{\pm e_1, \dots, \pm e_d\}$. In this case, let $j_n$ denote the number of susceptible neighbors of $v$ (including $a_n$), that is
\begin{equation*}
    j_n := \left\vert \{ u \in S_n : \langle v, u \rangle \in  E(\Z^d)  \} \right\vert .
\end{equation*}
Since we start with two infected vertices, necessarily $|S_n| \neq 2d$, and then $1 \leq j_n \leq 2d-1$. 


For each $i \in \{1, \dots, d\}$, we set the notation $\sigma(\pm e_i) = \pm u_i$.

Consider the event $A_n(v,e)$ where, either, for some $\ell \in \{1, \dots, 2d-1 \}$, the vertex $(x(v)+\sigma(e))^{(\ell)}$ is open or, for some $k \geq 1$ and $\ell \in \{1, \dots, 2d-1\}$,
\begin{itemize}
    \item the vertices $\left( x(v) + iu_{d+1} \right)^{(j_n)}$, $i = 1, \dots, k$, are open, and
    \item for each $i \in \{0, \dots, k-1\}$, the vertices $\left( x(v) + i u_{d+1} + \sigma(e) \right)^{(j)}$, $j = 1, \dots, 2d-1$, are closed, and
    \item the vertex $\left( x(v) + ku_{d+1} + \sigma(e) \right)^{(\ell)}$ is open. 
\end{itemize}
Observe that event $A_n$ is defined based on the whole previous construction up to time $n$, but is independent of the past. Therefore, the event $A_n(v,e)$ has probability   
\begin{eqnarray*}
    \Pp(A_n(v,e)) &=& \left(1 - (1-q)^{2d-1}\right) \sum_{i=0}^{\infty} \left[q(1-q)^{2d-1}\right]^i \\
    &=& \frac{p}{1 - (1-p)q}\\
    &=& \frac{p}{1-(1-p)\left(1-(1-p)^{1/(2d-1)}\right)} \\
    &=& \frac{p}{ p + (1-p)^{2d/(2d-1)} }.
\end{eqnarray*}

To conclude our induction step, we set 
\begin{itemize}
    \item $R_{n+1} = 
    \begin{cases} R_n \cup \{v+e\}, &\text{if } A_n(v,e) \text{ does not occur};\\
                  R_n, & \text{otherwise};
    \end{cases}$
    \item $I_{n+1} = 
    \begin{cases} I_n \cup \{v+e\}, &\text{if } A_n(v,e) \text{ occurs};\\
                  I_n, & \text{otherwise}.
    \end{cases}$
\end{itemize}
In the event $A_n(v,e)$, if for some $\ell \in \{1, \dots, 2d-1 \}$ the vertex $(x(v)+\sigma(e))^{(\ell)}$ is open, we set $x(v+e) =( x(v)+\sigma(e) )^{(\ell)}$. Otherwise, let $k \geq 1$ and $\ell \in \{1, \dots, 2d-1\}$ be such that the above three conditions are satisfied, in this case, we set $x(v+e) =( x(v)+ k u_{d+1} + \sigma(e) )^{(\ell)}$.

By construction, for each $v \in \cup_n I_n$, the projection of $x(v)$ onto $\Z^d$ is equal to $v$. Therefore $x: \Z^d \longrightarrow \Z^{d+1}_V$ is also injective. It follows that the site percolation model in $\Z^{d+1}_V$ (with $e_1^{(1)}$ declared open), dominates the infection process $(I_n)_{n \geq 0}$. 
Since  $|\cup_n I_n|$ has the same distribution as the size of the open cluster of $\{0, e_1\}$ in a supercritical percolation model on $\Z^d$ with parameter $p / [p + (1-p)^{2d/(2d-1)}] > p^s_c(d)$, the proof follows.
\end{proof}

\section{Proof of Theorems}\label{sec:proofs}
In this section, we prove the theorems using the couplings of the last section.

\subsection{Proof of Theorem \ref{theo:bounds}}

\begin{proof}
Note that in Proposition \ref{prop:coupling}, the selection of the final pair of coordinates was arbitrary; this approach can be applied to any pair of coordinates. Given $d \geq 3$, following this idea, we partition the set of directions of $\Z^d$ into three subsets of sizes $\left \lfloor\frac{d}{3} \right\rfloor, \left\lfloor\frac{d+1}{3}\right\rfloor, \left\lfloor\frac{d+2}{3}\right\rfloor$, and iteratively apply the proposition to each subset. One can readily verify that this partition provides the optimal bound. After $d-3$ applications of Proposition \ref{prop:coupling}, we obtain
\begin{equation*}
    \theta_d(p)\geq\theta_3(p_1,p_2,p_3),
\end{equation*}
where
\begin{equation*}
    p_i=1-(1-p)^{\lfloor \frac{d+i-1}{3}\rfloor} \: \text{ for } i=1,2,3.
\end{equation*} 
Now, by Proposition \ref{prop:triangular}, we have 
\begin{equation*}
    \theta_3(p_1,p_2,p_3)\geq \theta_{\mathbf{T}}(p_1,p_2,p_3).
\end{equation*}
Let $p^*(d)$ as defined as \eqref{eq:Theorem1}. From definition of $p^*(d)$, if $p>p^{\ast}(d)$ we have that $ p_1 + p_2 + p_3 - p_1p_2p_3 > 1$ and then, by Theorem 11.116 in \cite{Grim},  $\theta_{\mathbf{T}}(p_1,p_2,p_3)>0$.

In this way, we conclude that 
\begin{equation*} 
\theta_d(p)\geq\theta_3(p_1,p_2,p_3) \geq  \theta_{\mathbf{T}}(p_1,p_2,p_3)>0, 
\end{equation*}
which yields the result.
\end{proof}

\subsection{Proof of Theorem \ref{theo:bondO}}
\begin{proof}
First, we remark that  Proposition \ref{prop:coupling} holds {\it mutatis mutandis} for inhomogeneous oriented bond Bernoulli percolation. As a consequence, we have the following corollary which we state without proof.
\begin{lem}\label{lem:Cor}
Let $k \in \N$ be such that $d$ is divisible by $k$. Given $p \in [0,1]$, let $\tilde{p}$ be such that 
\[(1-\tilde{p}) = (1-p)^{d/k}.\]
Then,  $\vv{\theta}^{\, b}_d(p) \geq \vv{\theta}_k^{\, b}(\tilde{p})$.
\end{lem}

Taking $k=2$ in Lemma \ref{lem:Cor} and using  Liggett's upper bound $\vv{p}_c^{\, b}(2) \leq 2/3$ (see \cite{L}), the first item of Theorem \ref{theo:bondO} follows. The third item is equivalent to Proposition \ref{prop:crossover}. 

 Finally, the second item is implicitly proved in \cite{GPS}. 
 There, the authors defined a quantity $\lambda(1/d, \dots, 1/d)$ and showed that (see  Equation (3.2) in \cite{GPS}), if $p \in [0,1]$ satisfies
\[\phi(d) :=  \lambda(1/d, \dots, 1/d) \leq dp -1,\]
then $\vv{\theta}^{\, b}_d(p) > 0$. In particular, 
\[ \vv{p}_c^{\, b}(d) \leq \frac{1+\phi(d)}{d}.\]
The conclusion follows from the estimate given in Subsection 3.1 of \cite{GPS} (see the equation above (3.5)), that is
\[\phi(d) \leq \frac{1}{d}+\frac{8}{d^2}+\frac{d^{3/2}}{(\sqrt{2\pi})^{d-1}}\left[\frac{d-1}{d-3}\right]e^{\frac{1}{12d}}. \]
\end{proof}

\subsection{Proof of Theorem \ref{theo:siteNO}}
Inhomogeneous percolation is not well defined for site models, and therefore, we are not able to  formulate a version of Proposition \ref{prop:coupling} for site percolation. But note that Lemma \ref{lem:Cor} only involves homogeneous models. In the following, we give its site version.

\begin{lem}\label{lem:ksite}
Consider non-oriented site Bernoulli percolation on $\Z^{d}$. Let $k \in \N$ be such that $d$ is divisible by $k$. Given $p \in [0,1]$, let $\tilde{p}$ be such that 
\[(1-\tilde{p}) = (1-p)^{d/k}.\]
Then, we have ${\theta}^s_d(p) \geq {\theta}^s_k(\tilde{p})$.
\end{lem} 

\begin{proof}
The lemma follows by a coupling between non-oriented site Bernoulli percolation on $\Z^d$, with parameter $p$, and on $\Z^{k}$, with parameter $\tilde{p}$. To establish this coupling, we again construct a sequence of vertex sets $(I_n, x(I_n), R_n, S_n)_{n \geq 0}$.  

First, we recall that $\Ee = \{e_1, \dots, e_d\}$ is the set of positive unit vectors of $\Z^d$ and let $(D_{u_1}, \dots, D_{u_k})$ be a uniform partition of $\Ee$ into $k$ subsets indexed by the set $\{u_1, \dots, u_k \}\subset\Z^k$ of positive unit vectors of $\Z^k$. For each $u \in \{u_1, \dots, u_k \} $, let $D_{-u} = -D_u$.

Consider the non-oriented site Bernoulli percolation model on $\Z^d$, with parameter $p$. For each $v \in \Z^d$ and each $u \in \{\pm u_1, \dots, \pm u_k\}$, we define the following event
\begin{equation*}
    B(v,u) := \{ v+e \text{ is open, for some } e \in D_u \}.
\end{equation*}
If the event $B(v,u)$ occurs, we set $e(v,u) = v+e$, where $e$ is an open vertex that guarantee the occurrence of $B(v,u)$. Note that  $\Pp(B(v,u)) = \tilde{p}.$

For $n = 0$, we define 
\begin{itemize}
    \item $I_0 = \{0\} \subset \Z^k$;
    \item $R_0 = \emptyset \subset \Z^k$;
    \item $x(0) = 0 \in \Z^d$.
\end{itemize}
If for some $n \geq 0$, the sets $I_n, x(I_n)$ and $R_n$ are already defined, let $S_n$ be as given in \eqref{eq:SnSITE}. If $S_n = \emptyset$, the sequence becomes constant as in \eqref{eq:stop}. Otherwise, let $v_n$ be the smallest vertex of $S_n$ (according to a preset ordering of $\Z^k$). In this case, we write $v_n = v+u \notin I_n$, where $v \in I_n$ and $u \in \{\pm u_1, \dots, \pm u_k\}$. To conclude the induction step, we define
\begin{itemize}
    \item $R_{n+1} = 
    \begin{cases} R_n \cup \{v+u\}, &\text{if } B(x(v),u) \text{ does not occurs};\\
                  R_n, & \text{otherwise};
    \end{cases}$
    \item $I_{n+1} = 
    \begin{cases} I_n \cup \{v+u\}, &\text{if } B(x(v),u) \text{ occurs};\\
                  I_n, & \text{otherwise}.
    \end{cases}$
\end{itemize}
If $B(x(v),u)$ occurs, we define $x(v+u) = e(x(v),u)$.

Note  that, by construction, the function $x: \cup_n I_n \longrightarrow \Z^d$ is such that, for each $v = v_{u_1}u_1+ \cdots + v_{u_k}u_k \in \cup_n I_n$, writing $x(v) = x_{e_1} e_1 + \cdots + x_{e_d} e_d$, we have
\begin{equation*}
    \sum_{e \in D_u} x_e = v_u, \quad \forall u \in \{u_1, \dots, u_k\}.
\end{equation*}
Therefore, $x: \cup_n I_n \longrightarrow \Z^d$ is injective. The conclusion of the lemma follows as in previous couplings. 
\end{proof}

With Lemma \ref{lem:ksite} we are now able to conclude the goal of this section.

\begin{proof}[Proof of Theorem \ref{theo:siteNO}]
The first item of Theorem \ref{theo:siteNO} follows by taking $k=2$ in Lemma \ref{lem:ksite} toghether with Wierman's upper bound $p_c^s(2)\leq 0.68$ (see \cite{W2}). Taking $k=3$ in Lemma \ref{lem:ksite} and using the upper bound $p_c^s(3) < 1/2$ of Campanino-Russo (see \cite{CR}), the second item follows. The third item follows directly by Proposition \ref{prop:crossoversite}.
\end{proof}

\subsection{Proof of Theorem \ref{theo:siteO}}
\begin{proof}
First, we remark that Lemma \ref{lem:ksite} and Proposition \ref{prop:crossoversite} hold {\it mutatis mutandis} for oriented site Bernoulli percolation. Therefore to conclude the proof of Theorem \ref{theo:siteO}, the last ingredient is  Liggett's upper bound $\vv{p}^{\,s}_c(2) \leq 3/4$ (see \cite{L}).
\end{proof}


\acks The authors thank Roger~Silva for valuable comments on a first version of the manuscript and the two anonymous referees for their valuable comments and suggestions, which helped us improve the text.

\fund 
P.A.~Gomes has been supported by São Paulo Research Foundation (FAPESP), grant 2020/02636-3 and grant 2023/13453-5. A.~Pereira was partially funded by the Fundação de Amparo a Pesquisa de Alagoas,  Fapeal (Project APQ2022021000107; E:60030.0000000161/2022) and by Conselho Nacional de Desenvolvimento Científico e Tecnológico - CNPq, Universal (402952/2023-5). R.~Sanchis has been partially supported by Conselho Nacional de Desenvolvimento
Científico e Tecnológico (CNPq), CAPES and by FAPEMIG (APQ-00868-21 
 and RED-00133-21).

\competing 
There were no competing interests to declare which arose during the preparation or publication process of this article.



%
%
%

\end{document}